# Robust Unit Commitment Considering Strategic Wind Generation Curtailment

C. Wang, F. Liu, W. Wei, S. Mei
Department of Electrical Engineering
Tsinghua University
100084 Beijing, China
c-w12@mails.tsinghua.edu.cn

F. Qiu, J. Wang
Energy System Division
Argonne National Laboratory
60439 Argonne, IL, USA
fqiu@anl.gov

*Abstract*— Wind generation is traditionally treated as a non-dispatchable resource and is fully absorbed unless there are security issues. To tackle the operational reliability issues caused by the volatile and non-dispatchable wind generation, many dispatch frameworks have been proposed, including robust unit commitment (RUC) considering wind variation. One of the drawbacks that commonly exist in those dispatch frameworks is increased demand on flexibility resources and associated costs. To improve wind dispatchability and reduce flexibility resource costs, in this paper, we propose a novel RUC model considering strategic wind generation curtailment (WGC). Strategic WGC can reduce wind uncertainty and variability and increase the visibility of wind generation capacity. As a result, the ramping requirement for wind generation will be reduced and ramp-up capability of wind generation can be increased, leading to reduced day-ahead operational cost with guaranteed operational reliability requirement of power systems. The economic benefits also include profits gained by wind farm by providing ramping-up capacities other auxiliary services. We also propose a solution algorithm based on the column and constraint generation (C&CG). Simulations on the IEEE 39-bus system and two larger test systems demonstrate the effectiveness of the proposed RUC model and efficiency of the proposed computational methodology.

## NOMENCLATURE

*Indices*

| | |
|---|---|
| $g$ | Index for generators. |
| $m$ | Index for wind farms. |
| $l$ | Index for transmission lines. |
| $j$ | Index for loads. |
| $t$ | Index for time periods. |

*Parameters*

| | |
|---|---|
| $T$ | Number of time periods. |
| $M$ | Number of wind farms. |
| $G$ | Number of thermal generators. |
| $L$ | Number of transmission lines. |
| $S_g$ | Start-up cost of generator $g$. |
| $c_g$ | Constant term of generation cost function. |
| $P_g^{min}/P_g^{max}$ | Minimal/ maximal output of generator $g$. |
| $R_+^g/R_-^g$ | Ramp-up/ ramp-down limit for generator $g$. |
| $T_g^{on}/T_g^{off}$ | Minimum on/off hour of generator $g$. |
| $F_l$ | Transmission capacity of line $l$. |
| $W$ | Wind generation uncertainty set. |
| $\hat{w}_{mt}$ | Forecasted output of wind farm $m$ in period $t$. |
| $\Gamma^S/\Gamma^T$ | Uncertainty budget over spatial/ temporal scale. |
| $w_{mt}^u/w_{mt}^l$ | Upper/ lower bound of $W$. |
| $D_{jt}$ | Load demand of load node $j$ in period $t$. |
| $\beta_t/\beta_s$ | Confidence level of $\Gamma^T/\Gamma^S$. |
| $\pi_{gt}/\pi_{mt}/\pi_{jt}$ | Generation shift distribution factor of generator $g$/ wind farm $m$/ load $j$ in period $t$. |

*Decision Variables*

| | |
|---|---|
| $\alpha_{mt}$ | Wind generation commitment ratio. |
| $u_{gt}$ | Binary variable indicating whether generator $g$ is on or off in period $t$. |
| $z_{gt}$ | Binary variable indicating whether generator $g$ is started up in period $t$. |
| $P_{gt}$ | Real-time output of generator $g$ in period $t$. |
| $\hat{P}_{gt}$ | Day-ahead output of generator $g$ in period $t$. |
| $v_{mt}^u/v_{mt}^l$ | Binary variable indicating normalized positive /negative output deviation. |
| $\Delta w_{mt}$ | Wind generation curtailment in the recourse stage. |
| $\Delta D_{jt}$ | Load shedding in the recourse stage. |

## I. INTRODUCTION

The increasing penetration of wind generation has brought many challenges to power system, especially in the scope of operation. On one hand, wind generation cannot be forecasted accurately. According to [1], the day-ahead wind generation forecast error can be 20% or larger. On the other hand, wind generation is traditionally treated as non-dispatchable, which means it is fully absorbed unless there are security issues. To cope with these tough characters of wind generation, operation strategies with more flexibility are desiderated. One of the focuses is Unit commitment (UC) as it determines the operational flexibility of the following day.

Robust unit commitment (RUC) has been an active topic recently as it can guarantee the operational security given a prescribed wind generation uncertainty set, which is determined by the operational reliability requirements (ORR) of power

This work is supported in part by the Foundation for Innovative Research Groups of the NNSF of China (51321005), and State Key Development Program of Basic Research of China (2013CB228201).

systems [2], [3]. To meet ORR, expensive fast-response generators are dispatched more frequently under RUC to provide sufficient regulation capability, which significantly increases the operational cost of RUC. Numerous models and approaches have been proposed to reduce the operational cost of RUC, such as minimax regret RUC [4], unified stochastic and robust UC [5], and hybrid stochastic/ interval RUC [6]. However, these approaches mainly focus on exploring operational flexibility from conventional generators.

As a matter of fact, wind farms can contribute to operational flexibility of power systems. In [7], wind farms are treated as a reserve provider by decreasing its output and the system dispatch cost is decreased because of the reduction of reserve requirement. Wind farms can also benefit from this practice as its uncommitted power can be used to provide ramping-up reserve. The incentives of wind farms of providing reserve service in the perspective of electricity market are analyzed in [8], [9]. Indeed, decreasing the output of wind farms, also called wind generation curtailment (WGC), is a useful method to recover operational feasibility.

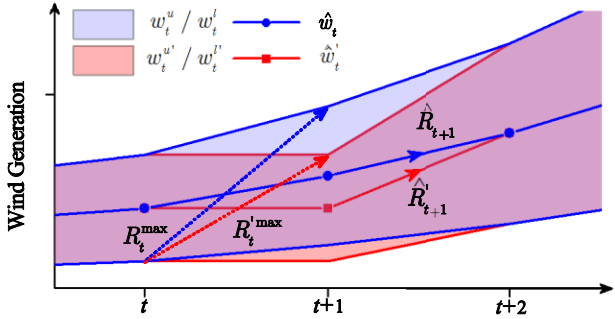

Fig. 1. Schematic diagram of RUC with strategic WGC.

In this paper, a novel RUC formulation considering strategic WGC is proposed, aiming to reduce the day-ahead dispatch costs. The schematic diagram of strategic WGC is shown in Fig. 1. The blue lines represent forecasted value of wind generation $\hat{w}_t$ and the red lines represent the committed wind generation $\hat{w}'_t$. The blue and red dash area represent the boundaries of uncertainty set of $\hat{w}_t$ and $\hat{w}'_t$, respectively, according to operational reliability requirement of power systems. Intuitively, strategic WGC occurs in period $t$+1. This practice may benefit the power system from Four aspects: (1) Reducing ramping requirement. After strategic WGC, the maximum ramping requirement of wind generation between period $t$ and $t$+1 decreases from $R_t^{max}$ to $R_t'^{max}$. (2) Increasing ramping up capability. The ramping up capacity provided by wind generation between period $t$+1 and $t$+2 increases from $\hat{R}_{t+1}$ to $\hat{R}'_{t+1}$. (3) Improving utilization ratio of transmission lines. As the uncertainty band becomes narrower compared with original one in period $t$+1, transmission lines can leave less space for possible over-generation of wind generation. (4) Enhancing solvability of RUC. In traditional RUC, robust dispatch strategy may not be generated, especially under critical operational reliability requirement. With strategic WGC, the solvability of RUC can always be guaranteed.

Meanwhile, wind farms can also benefit from strategic WGC: on one hand, it can charge power system for providing ramping up capacity; on the other hand, it can use the uncommitted or the "curtailed" wind generation to provide auxiliary services. Detailed discussion can be found in Section IV. C.

The rest of the paper is structured as follows: Section II describes the mathematical formulation of RUC with strategic WGC. Section III presents the solution methodology. In Section IV, illustrative examples on different test systems are proceeded to demonstrate the effectiveness of the proposed model and algorithm. Finally, Section V gives the conclusion of the paper.

## II. MATHEMATICAL FORMULATION

### A. RUC Considering Strategic WGC

$$\min_{z,u,\hat{p},\alpha} \sum_{t=1}^{T}\sum_{g=1}^{G}\left(S_g z_{gt} + c_g u_{gt} + C_g(\hat{p}_{gt})\right) \quad (1a)$$

$$s.t.\ -u_{g(t-1)} + u_{gt} - u_{gk} \le 0,\ \forall g, \forall t, k=t,\cdots,t+T_g^{on}-1. \quad (1b)$$

$$u_{g(t-1)} - u_{gt} + u_{gk} \le 1,\ \forall g, \forall t, k=t,\cdots,t+T_g^{off}-1. \quad (1c)$$

$$-u_{g(t-1)} + u_{gt} - z_{gt} \le 0,\ \forall g, \forall t \quad (1d)$$

$$u_{gt}P_{\min}^g \le \hat{p}_{gt} \le u_{gt}P_{\max}^g\ \forall g, \forall t \quad (1e)$$

$$\hat{p}_{gt} - \hat{p}_{g(t+1)} \le u_{g(t+1)}R_-^g + (1-u_{g(t+1)})P_{\max}^g\ \forall g,\ \forall t \quad (1f)$$

$$\hat{p}_{g(t+1)} - \hat{p}_{gt} \le u_{gt}R_+^g + (1-u_{gt})P_{\max}^g\ \forall g,\ \forall t \quad (1g)$$

$$\sum_{g=1}^G \hat{p}_{gt} + \sum_{m=1}^M \alpha_{mt}\hat{w}_{mt} = \sum_{j=1}^J D_{jt}\ \forall t \quad (1h)$$

$$-F_l \le \sum_{g=1}^G \pi_{gt}\hat{p}_{gt} + \sum_{m=1}^M \pi_{mt}\alpha_{mt}\hat{w}_{mt} - \sum_{j=1}^J \pi_{jt}D_{jt} \le F_l\ \forall t \quad (1i)$$

$$0 \le \alpha_{mt} \le 1\ \forall m\ \forall t \quad (1j)$$

$$u_{gt}, \alpha_{mt} \in \Omega \quad (1k)$$

In model (1), (1a) is to minimize the day-ahead operational cost, in which the first term represents the start-up cost of units and the last two terms represent the economic dispatch (ED) cost under base case. In (1a), $C_g(\cdot)$ is quadratic and can be further linearized by piecewise linearization method. (1b) and (1c) describes the minimum on/off period limits of generators. (1d) is the start-up constraints of generators. (1e) is the generation capacity of generators. (1f) and (1g) are the ramping rate limits of generators, respectively. (1h) depicts the power balance requirement under base case. (1i) is the network power flow limits on transmission lines. (1j) depicts the upper and lower boundary of $\alpha_{mt}$. $\Omega$ is the feasibility set of $u_{gt}, \alpha_{mt}$ and its definition is as follows.

$$\Omega:=\left\{\max_{v^u,v^l}\min_{p,\Delta w,\Delta D}\sum_{t=1}^T\left(\sum_{m=1}^M\Delta w_{mt}+\sum_{j=1}^J\Delta D_{jt}\right)=0\right. \quad (2a)$$

$$s.t.\ u_{gt}P_{\min}^g \le p_{gt} \le u_{gt}P_{\max}^g\ \forall g, \forall t \quad (2b)$$

$$p_{gt} - p_{g(t+1)} \le u_{g(t+1)}R_-^g + (1-u_{g(t+1)})P_{\max}^g\ \forall g, \forall t \quad (2c)$$

$$p_{g(t+1)} - p_{gt} \le u_{gt}R_+^g + (1-u_{gt})P_{\max}^g\ \forall g, \forall t \quad (2d)$$

$$\sum_{g=1}^G p_{gt} + \sum_{m=1}^M (w_{mt}-\Delta w_{mt}) = \sum_{j=1}^J (D_{jt}-\Delta D_{jt}) \quad (2e)$$

$$0 \le \Delta D_{jt} \le D_{jt}\ \forall j, \forall t \quad (2f)$$

$$0 \le \Delta w_{mt} \le w_{mt}\ \forall m, \forall t \quad (2g)$$

$$-F_l \leq \sum_{g=1}^{G} \pi_{gt} p_{gt} + \sum_{m=1}^{M} \pi_{mt}(w_{mt} - \Delta w_{mt}) - \cdots$$
$$- \sum_{j=1}^{J} \pi_{jt}(D_{jt} - \Delta D_{jt}) \leq F_l \quad \forall t \quad (2h)$$

$$w_{mt} = (w_{mt}^u - \hat{w}_{mt})\alpha_{mt} v_{mt}^u + (w_{mt}^l - \hat{w}_{mt})\alpha_{mt} v_{mt}^l + \alpha_{mt}\hat{w}_{mt} \quad (2i)$$

$$\sum_{t=1}^{T}(v_{mt}^u + v_{mt}^l) \leq \Gamma^T \quad \forall m \quad (2j)$$

$$\sum_{m=1}^{M}(v_{mt}^u + v_{mt}^l) \leq \Gamma^S \quad \forall t \quad (2k)$$

$$v_{mt}^u + v_{mt}^l \leq 1 \quad \forall m, \forall t \quad (2l)$$

$$v_{mt}^u, v_{mt}^l \in \{0,1\} \quad (2m)$$

In model (2), (2a) is the sum of load shedding (LS) and WGC. (2b) depicts the capacity of generators. (2c) and (2d) limit the ramping capacity of generators. (2e) depicts the relaxed power balance requirement with emergency actions including LS and WGC. (2f) and (2g) are the boundary of LS and WGC respectively. (2h) is the network power flow limits considering LS and WGC. (2i)-(2m) use a polyhedral set to describe the wind generation denoted as *W*. Specifically, (2i) depicts the wind generation output; (2j) and (2k) describe the uncertainty budgets over both temporal and spatial domains, respectively.

It should be pointed out that WGC in model (2) occurs in the recourse stage, which strategic WGC is a day-ahead decision. From model (1) and (2), the proposed RUC is a two-stage robust optimization problem. The first stage decision variables are $z_{gt}, u_{gt}, \hat{p}_{gt}, \alpha_{mt}$, the recourse action variables are $p_{gt}, \Delta w_{mt}, \Delta D_{jt}$, and the uncertainty variables are $v_{mt}^u, v_{mt}^l$. Due to the existence of (2a), there will no LS nor WGC within the uncertainty set *W*, which guarantees the operational feasibility of $u_{gt}$ and $\alpha_{mt}$.

### B. Compact Model

For simplicity, the compact formulation of the proposed RUC with strategic WGC can be written as follows:

$$\min_{\mathbf{x},\hat{\mathbf{y}},\boldsymbol{\alpha}} \mathbf{b}^T\mathbf{x} + \mathbf{c}^T\hat{\mathbf{y}} \quad (3a)$$

$$s.t. \quad \mathbf{Bx} + \mathbf{C}\hat{\mathbf{y}} + \mathbf{D}\boldsymbol{\alpha} \leq \mathbf{d} \quad (3b)$$

$$\begin{pmatrix}\mathbf{x}\\ \boldsymbol{\alpha}\end{pmatrix} \in \begin{cases} \max_{\mathbf{v}} \min_{\mathbf{y},\mathbf{s}} \mathbf{f}^T\mathbf{s} = 0 & (3c)\\ s.t. \quad \mathbf{Ex} + \mathbf{Fy} + \mathbf{G}(\boldsymbol{\alpha}\circ\mathbf{v}) + \mathbf{Hs} + \mathbf{Jv} \leq \mathbf{g} & (3d)\\ \mathbf{Lv} \leq \mathbf{h} & (3e)\end{cases}$$

In model (3), **x** represents the binary vector of generators. $\hat{\mathbf{y}}, \mathbf{y}$ represent the vector of generators. **α** represents wind generation commitment ratio vector. **s** represents the LS and WGC vector. **v** is the binary vector depicting wind generation uncertainty. **b, c, d, f, g, h, B, C, D, E, F, G, H, J, L** are constant coefficient matrix and can be derived from model (1) and (2). Specially, $\boldsymbol{\alpha} \circ \mathbf{v}$ is a Hadamard product.

### III. SOLUTION METHODLOGY

In this section, we will derive the solution methodology to solve problem (3), in which (3a)-(3b) depicts the main problem (MP) and (3c)-(3e) formulates the feasibility-checking subproblem.

### A. Solution Methodology for Subproblem

The subproblem is a bi-level mixed integer linear problem (MILP) and can be solved by many effective methods, such as the Karush-Kuhn-Tucker (KKT) condition based method [10], and the strong duality theory based method [11]. In this paper, the inner problem of (3c) is replaced by its dual problem to reformulate (3c) as a single-level bilinear program, which can be furthered solved by big-M linearization method [12]. As the big-M linearization method is proved effective with high efficiency and accuracy in practice, this paper adopts it to solve (3c)-(3e). The compact formulation of dual problem of (3c) is as follows.

$$\max_{\mathbf{v},\boldsymbol{\lambda}} R = \boldsymbol{\lambda}^T(\mathbf{g}-\mathbf{Ex}) - \boldsymbol{\lambda}^T\mathbf{Jv} - \boldsymbol{\lambda}^T\mathbf{G}(\boldsymbol{\alpha}\circ\mathbf{v}) \quad (4a)$$

$$s.t. \quad [\mathbf{F} \vdots \mathbf{H}]^T \boldsymbol{\lambda} \leq [\mathbf{0}^T \vdots \mathbf{f}^T]^T \quad (4b) \quad (4)$$

$$\boldsymbol{\lambda} \leq \mathbf{0} \quad (4c)$$

$$(3e)$$

Where, **λ** is the dual variable vector of inner problem of (3c). Noticed that there are bilinear terms in (4a), auxiliary variables and constraints are introduced to replace them and (4) can be transferred into a MILP problem as follows.

$$\max_{\mathbf{v},\boldsymbol{\lambda},\boldsymbol{\gamma}} R = \boldsymbol{\lambda}^T(\mathbf{g}-\mathbf{Ex}) - \boldsymbol{\gamma}^T\mathbf{q} \quad (5a)$$

$$s.t. \quad (3e), (4b)-(4c)$$

$$-M_{Big}\mathbf{v} \leq \boldsymbol{\gamma} \leq \mathbf{0} \quad (5b) \quad (5)$$

$$-M_{Big}(\mathbf{1}-\mathbf{v}) \leq \boldsymbol{\lambda} - \boldsymbol{\gamma} \leq \mathbf{0} \quad (5c)$$

Where, **γ** is the auxiliary variable vector, **q** is a constant vector and can be derived from the following formula.

$$\boldsymbol{\lambda}^T\mathbf{Jv} + \boldsymbol{\lambda}^T\mathbf{G}(\boldsymbol{\alpha}\circ\mathbf{v}) = \sum_i\sum_j q_{ij}\lambda_i v_j = \boldsymbol{\gamma}^T\mathbf{q}, \quad \gamma_{ij} = \lambda_i v_j \quad (6)$$

(5b) and (5c) are auxiliary constraints generated during objective function linearization using the big-M method. $M_{big}$ is sufficient large positive real number. Finally, the subproblem can be transformed into a single-level MILP with additional binary variables and constraints with big-M parameter, which can be solved easily by commercial solvers such as CPLEX.

### B. Solution Methodology for Main Problem

Noted that both (3a)-(3b) and the subproblem are MILPs. Next a column and constraint generation (C&CG) based algorithm is developed to solve (3) and named as A1. The details of A1 is as follows.

| A1: C&CG-based Algorithm |
|---|
| *Step 1:* set $i=0$ and $\mathbf{O} = \emptyset$. |
| *Step 2:* Solve (3a)-(3b) with the additional constraints as follows. |
| $\quad\quad \mathbf{Ex} + \mathbf{Fy}^k + \mathbf{G}(\boldsymbol{\alpha}\circ\mathbf{v}_k^*) + \mathbf{Jv}_k^* \leq \mathbf{g} \quad \forall k \leq i \quad (7a)$ |
| *Step 3:* Solve the subproblem. If $|R_{k+1} - R_k| < \epsilon$, terminate. Otherwise, derive the optimal solution $\mathbf{v}_{k+1}^*$, create variable vector $\mathbf{y}^{k+1}$ and add the following constraints |
| $\quad\quad \mathbf{Ex} + \mathbf{Fy}^{k+1} + \mathbf{G}(\boldsymbol{\alpha}\circ\mathbf{v}_{k+1}^*) + \mathbf{Jv}_{k+1}^* \leq \mathbf{g} \quad (7b)$ |
| *Update* $i=i+1$, $\mathbf{O} = \mathbf{O}\cup\{i+1\}$ *and go to Step 2.* |

In A1, $\epsilon$ represents the convergence gap and $R_k$ represents the optimal value of subproblem in *k*th iteration. In traditional

C&CG algorithm [13], a set of constraints (3d) of subproblem with the identified worst-case scenario are directly added to MP. However, in A1, the added constraints (7b) are not the same with the original constraints (3d) in subproblem. Compared with (7b), (3d) can be regarded as loose constraints with slack variables as emergency regulation is involved. This difference makes A1 a C&CG-based algorithm.

## IV. ILLUSTRATIVE EXAMPLE

In this section, we present numerical experiments carried on the modified IEEE 39-bus system, the modified IEEE 118-bus system and a regional power grid of China to show the effectiveness of the proposed model and algorithms. The experiments are performed on a PC with Intel(R) Core(TM) 2 Duo 2.2 GHz CPU and 4 GB memory. The optimality gap is set as 0.1% in this section.

### C. The Modified IEEE 39-bus System

The tested system has 10 generators and 46 transmission lines. A wind farm is connected to the system at bus 29 with an installed capacity of 500 MW. The generators' parameters can be found in [14]. The load curve and the day-ahead forecast of wind generation are both scaled down from the day-ahead curve of California ISO as shown in Fig. 2. We choose the confidence level $\beta_t = 95\%$, yielding $\Gamma^T \approx 8$ [17]. The standard deviation of wind generation forecast error is subject to (8) with $\sigma_m = 0.15$ and its mean value is zero. In this case, wind generation within 99% forecast error band is required to be fully accommodated and the forecast error bands are derived by Gaussian distribution for simplicity.

$$\sigma_{mt} = \sigma_m \cdot \hat{w}_{mt} \cdot (1 + e^{-(T-t)}) \quad \forall m, \forall t \qquad (8)$$

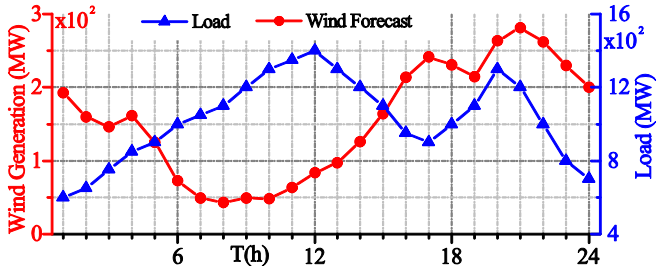

Fig. 2. Data of load demand and forecasted wind power.

### D. Comparision with Traditional RUC

In this subsection, RUC with strategic WGC is compared with traditional RUC in terms of operational cost under different penetration levels of wind generation. Specifically, the tradition RUC model is from [3]. Operational cost under different RUC models as well as wind penetration levels are listed in Table I. From Table I, the operational cost of RUC with strategic WGC is always lower than traditional RUC under the same wind penetration level and the minimum operational cost gap between two RUC models is around 1.4%. Meanwhile, as wind penetration level increases, the operational cost of traditional RUC decreases in the beginning and increases rapidly afterwards. Particularly, there will be no robust dispatch strategy when wind generation reaches 150% under traditional RUC. The reason is as follows, higher wind penetration level results in wider uncertainty band as the reliability requirement remains the same, which requires more critical ramping capability of power systems. However, the operational cost of RUC with strategic WGC keeps decreasing as wind generation penetration level increases, which reflects the superiority of the proposed RUC model compared to traditional ones, especially under high wind penetration case.

The optimal wind generation commitment ratio (WGCR) and the committed wind generation under 100% wind generation penetration is shown in Fig. 3. From Fig. 3, the committed wind generation (WG) is less than the forecasted WG in some periods, i.e., period 3, 4, 16, 17, 21, 23 and 24. Also, the corresponding error band (EB) in those periods are narrower than the original ones, which equals to less ramping capability requirement and contributes to the operational cost decrement.

TABLE I. OPERATIONAL COST UNDER DIFFERENT RUC MODELS

|  | Penetration Level (%) | Total ($) | UC ($) | ED ($) |
|---|---|---|---|---|
| Traditional RUC | 100 | 4.476×10⁵ | 1.33×10⁴ | 4.343×10⁵ |
|  | 110 | 4.449×10⁵ | 1.50×10⁴ | 4.299×10⁵ |
|  | 120 | 4.489×10⁵ | 2.50×10⁴ | 4.239×10⁵ |
|  | 130 | 4.549×10⁵ | 3.45×10⁴ | 4.205×10⁵ |
|  | 140 | 4.607×10⁵ | 4.13×10⁴ | 4.194×10⁵ |
|  | 150 | No Solution | - | - |
| RUC with Strategic WGC | 100 | 4.415×10⁵ | 7.66×10³ | 4.339×10⁵ |
|  | 110 | 4.381×10⁵ | 7.60×10³ | 4.305×10⁵ |
|  | 120 | 4.346×10⁵ | 7.66×10³ | 4.269×10⁵ |
|  | 130 | 4.322×10⁵ | 7.14×10³ | 4.250×10⁵ |
|  | 140 | 4.295×10⁵ | 8.12×10³ | 4.214×10⁵ |
|  | 150 | 4.275×10⁵ | 8.82×10³ | 4.187×10⁵ |

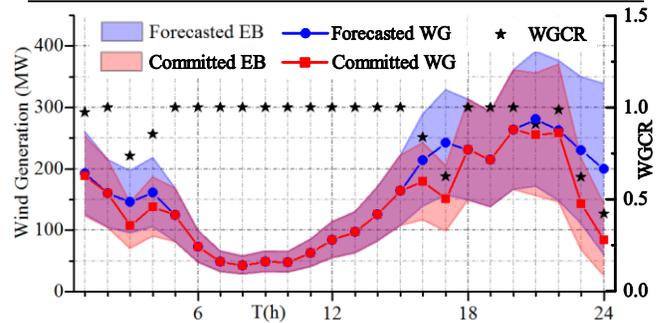

Fig. 3. Committed wind generation under 100% wind generation penetration.

### E. Incentives for Wind Farms

In this subsection, the reasons why wind farms are willing to accept WCGR are analyzed. One reason is that wind farms can benefit more from the committed wind generation by providing ramping product, as shown in Fig. 4. From Fig. 4, the ramp capability of RUC with strategic WGC is superior to traditional RUC in some periods, i.e., period 2, 3, 17, 23 and 24. The average ramping up and down capability of RUC with strategic WGC are 13.0 MW/h and 17.5 MW/h, respectively; while those of traditional RUC are 12.2 MW/h and 10.6 MW/h, respectively. Another reason, wind farms can make full use of the uncommitted wind generation by providing auxiliary services such as upward reserves and the expectation benefit $Q_m$ subject to (9), where $\delta_{mt}$ represents the price coefficient. Further, system operators can design some mechanism to compensate wind farms in order to encourage them accepting strategic WGC.

$$Q_m = \sum_{t=1}^{T} \delta_{mt}(1-\alpha_{mt})\hat{w}_{mt} \quad \forall m \tag{9}$$

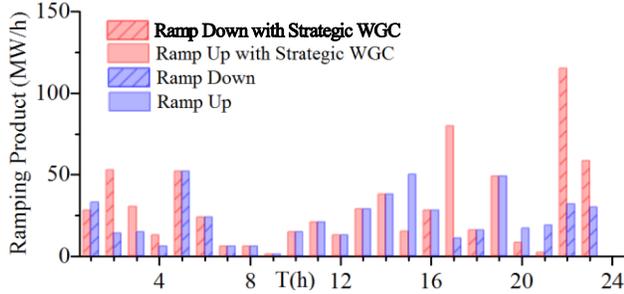

Fig. 4. Ramping product providing capability under different RUC models.

### F. Application Perspective

In this subsection, computational efficiency of the proposed model on larger test systems will be analyzed, as listed in Table II. The data of the modified IEEE 118-bus system can be found in [15], in which three wind farms are integrated; the data of Guangdong power grid can be found in [16], in which four wind farms are integrated. In Table II, the computational efficiency of traditional RUC model are listed as benchmark. From Table II, the proposed model will not increase the computational burden compared with traditional RUC model. In contrast, it can decrease the iteration number as well as computational time as the dispatchbility of wind farms are enhanced via introducing WGCR.

TABLE II. COMPUTATIONAL EFFICIENCY UNDER LARGER TEST SYSTEMS

|  |  | Iteration | Time (s) |
|---|---|---|---|
| IEEE 118-bus System | Traditional RUC | 4 | 137 |
|  | RUC with Strategic WGC | 3 | 86 |
| Guangdong Power Grid | Traditional RUC | 3 | 1317 |
|  | RUC with Strategic WGC | 3 | 1109 |

## V. CONCLUSION

In this paper, a RUC model considering strategic WGC is proposed to cope with the volatile and non-dispatchable wind generation. WGCR is introduced to control the proportion of committed wind generation. Compared with traditional RUC models, the proposed RUC with strategic WGC has three advantages. Economically, strategic WGC would benefit power systems by reducing ramping requirement, increasing ramping up capability as well as improving the utilization ratio of transmission lines, all of which would contribute to operational cost decrement of power systems. From solvability aspect, traditional RUC may fail to generate a robust dispatch strategy under critical operational reliability requirement, however, the solvability of RUC with strategic WGC could always be guaranteed. In terms of computational efficiency, the proposed model somewhat increases the dispatchbility of the non-dispatchable wind generation by introducing WGCR to the first stage of RUC, which results in better computational efficiency than traditional RUC. Mathematically, the RUC with strategic WGC is still a two-stage robust optimization problem, which can be efficiently solved by the proposed C&CG-based iterative algorithm. Simulations are carried out on the modified IEEE 39-bus system to illustrate the effectiveness of the proposed model and algorithm. The proposed methodology is also applied to IEEE 118-bus system and the real Guangdong Power Grid, demonstrating the practicality of our methodology.

## REFERENCES

[1] G. Giebel, R. Brownsword, G.Kariniotakis, M. Denhard, and C. Draxl, The State of the Art in Short-Term Prediction of Wind Power: A Literature Overview, 2nd edition [Online]. Available: http://www.anemosplus.eu/

[2] D. Bertsimas, E. Litvinov, and X. Sun, "Adaptive Robust Optimization for the Security Constrained Unit Commitment Problem," *IEEE Trans. Power Syst.*, vol.28, no.1, pp.52-63, Feb. 2013.

[3] R. Jiang, J. Wang, and Y. Guan, "Robust Unit Commitment with Wind Power and Pumped Storage Hydro," *IEEE Trans. Power Syst.*, vol.27, no.2, pp.800-810, May 2012.

[4] R. Jiang, J. Wang, and M. Zhang, "Two-Stage Minimax Regret Robust Unit Commitment," *IEEE Trans. Power Syst.*, vol.28, no.3, pp.2271-2282, Aug. 2013.

[5] C. Zhao, and Y. Guan, "Unified Stochastic and Robust Unit Commitment," *IEEE Trans. Power Syst.*, vol.28, no.3, pp.3353-3361, Aug. 2013.

[6] Y. Dvorkin, H. Pandzic, and M. A. Ortega-Vazquez, "A Hybrid Stochastic/Interval Approach to Transmission-Constrained Unit Commitment," *IEEE Trans. Power Syst.*, vol.30, no.2, pp.621-631, Mar. 2015.

[7] Y. Dvorkin, M. A. Ortega-Vazquez, and D. S. Kirschen, "Wind generation as a reserve provider," *IET Gener. Transm. Dis.*, vol.9, no.8, pp. 779-787, Nov. 2014.

[8] J. Liang, S. Grijalva, and R. G. Harley, "Increased Wind Revenue and System Security by Trading Wind Power in Energy and Regulation Reserve Markets," *IEEE Trans. Sust. Ener.,* vol. 2, no.3, pp. 340-347, Jul. 2011.

[9] T. Soares, P. Pinson, and H. Morais, "Optimal Offering Strategies for Wind Power in Energy and Primary Reserve Markets," available: "http:// pierrepinson.com/docs/soaresetal2015.pdf".

[10] S. J. Kazempour, A. J. Conejo, and C. Ruiz, "Strategic Generation Investment Using a Complementarity Approach," *IEEE Trans. Power Syst.*, vol.26, no.2, pp.940-948, May 2011.

[11] J. M. Arroyo, "Bilevel programming applied to power system vulnerability analysis under multiple contingencies," *IET Gener. Transm. Dis.*, vol.4, no.2, pp.178-190, Feb. 2010.

[12] B. Zeng, and L. Zhao, "Robust unit commitment problem with demand response and wind energy," in *Proc. IEEE Power Energy Society General Meeting*, San Diego, CA, USA, Jul. 2012.

[13] B. Zeng, and L. Zhao, "Solving two-stage robust optimization problems using a column-and-constraint generation method," *Oper. Res. Lett.*, vol.41, no.5, pp.457 - 461, 2013.

[14] Available: http://sys.elec.kitami-it.ac.jp/ueda/demo/WebPF/39-New-England.pdf.

[15] Available: "http://motor.ece.iit.edu/data/JEAS_IEEE118.doc".

[16] W. Wei, F. Liu, and S. Mei, "Two-level unit commitment and reserve level adjustment considering large-scale wind power integration," *Int. Trans. Electr. Energ. Syst.*, vol.24, no.12, pp.1726-1746, Oct. 2013.